\documentclass[11pt]{amsart}

\title{The Positive Bergman Complex of an Oriented Matroid}

\author{Federico Ardila, Caroline Klivans, and Lauren Williams}

\address{Mathematical Sciences Research Institute, Berkeley, CA 94720,
Computing and Information Science, Cornell University, Ithaca, NY 14853,
Department of Mathematics, MIT, Cambridge, MA 02139}

\usepackage{pstricks, pst-node, amssymb}
\psset{unit=1pt, arrowsize=4pt, linewidth=.7pt}
\psset{linecolor=blue}
\newgray{grayish}{.90}
\newrgbcolor{embgreen}{0 .5 0}

\def\vblack(#1, #2)#3{\cnode*[linecolor=black](#1, #2){3}{#3}}
\def\vwhite(#1,#2)#3{\cnode[linecolor=black,fillcolor=white,fillstyle=solid](#1,#2){3}{#3}}
\countdef\x=23
\countdef\y=24
\countdef\z=25
\countdef\t=26

\def\tbox(#1,#2)#3{
\x=#1 \y=#2
\multiply\x by 12
\multiply\y by 12
\z=\x \t=\y
\advance\z by 12
\advance\t by 12
\psline(\x,\y)(\x,\t)(\z,\t)(\z,\y)(\x,\y)
\advance\x by 6
\advance\y by 6
\rput(\x,\y){{\bf #3}}}

\usepackage{amsmath, amsthm, amssymb, amsbsy}
\usepackage{amsfonts, latexsym, stmaryrd, amscd, xy}
\usepackage[mathscr]{eucal}
\usepackage{epsfig}
\newtheorem{theorem}{Theorem}[section]

\newtheorem{proposition}[theorem]{Proposition}

\newtheorem{example}[theorem]{Example}
\newtheorem{corollary}[theorem]{Corollary}

\usepackage{amsfonts}
\usepackage{xy}
\usepackage{amssymb}
\usepackage{amsmath}
\newcommand{\R}{\mathbb{R}}
\newcommand{\w}{\omega}
\newcommand{\B}{{\mathcal{B}}}

\newcommand{\RR}{{\mathcal{R}}}

\newcommand{\FF}{{\mathcal{F}}}

\DeclareMathOperator{\Trop}{Trop}

\DeclareMathOperator{\init}{in} 

\newcommand{\thmrefer}[1]{\renewcommand\thetheorem
  {\protect\ref{#1}}\addtocounter{theorem}{-1}}

\xyoption{all} \CompileMatrices

\begin{document}

\begin{abstract}
We study the positive Bergman complex $\B^+(M)$
of an oriented matroid $M$, which is a certain subcomplex of the Bergman
complex $\B(\underline{M})$
of the underlying unoriented matroid $\underline{M}$.  The positive
Bergman complex is defined so that given a linear ideal $I$ with
associated oriented matroid $M_I$, the positive tropical variety
associated to $I$ is equal to the fan over ${\B}^+(M_I)$.  Our
main result is that a certain ``fine" subdivision of $\B^+(M)$ is
a geometric realization of the order complex of the proper part of
the Las Vergnas face lattice of $M$.  It follows that $\B^+(M)$ is
homeomorphic to a sphere.  For the oriented matroid of the
complete graph $K_n$, we show that the face poset of the ``coarse"
subdivision of $\B^+(K_n)$ is dual to the face poset of the
associahedron $A_{n-2}$, and we give a 
formula for the number
of fine cells within a coarse cell.
\end{abstract}

\maketitle

\section{Introduction}

In~\cite{Bergman}, Bergman defined the \emph{logarithmic
limit-set} of an algebraic variety in order to study its
exponential behavior at infinity. We follow \cite{Sturmfels} in
calling this set the \emph{Bergman complex} of the variety.
Bergman complexes have recently received considerable attention in
several areas, such as tropical algebraic geometry and dynamical
systems. They are the \emph{non-Archimedean amoebas} of
\cite{Einsiedler} and the \emph{tropical varieties} of \cite{
Speyer, Sturmfels}.

When the variety is a linear space, so that the defining ideal $I$
is generated by linear forms, Sturmfels \cite{Sturmfels} showed
that the Bergman complex can be described solely in terms of the
matroid associated to the linear ideal.  He used this description
to define the Bergman complex $\B(M)$ of an arbitrary matroid $M$.
Ardila and Klivans ~\cite{CarlyFed} showed that, appropriately subdivided,
the Bergman complex of a matroid $M$ is the order complex of the
proper part of the lattice of flats $L_M$ of the matroid.  This
result implies that the Bergman complex of an arbitrary matroid
$M$ is a finite, pure polyhedral complex, which is homotopy
equivalent to a wedge of spheres.

Total positivity is another topic which has received a great deal
of recent interest.  Although the classical theory concerns
matrices in which all minors are positive, in the past decade this
theory has been extended by Lusztig ~\cite{Lusztig1, Lusztig2},
who introduced the totally positive variety
$G_{>0}$ in an arbitrary reductive group $G$ and the totally
positive part $B_{>0}$ of a real flag variety $B$.  More recently,
the positive part of the tropicalization of an affine variety (or
positive tropical variety, for short) was
introduced by Speyer and Williams ~\cite{SpeyerWilliams}.

Sturmfels~\cite{Sturmfels2} suggested the notion of a positive
Bergman complex $\B^+(M)$ of an oriented matroid $M$ and
conjectured its relation to the Las Vergnas face lattice of $M$.
We define the positive Bergman complex and positive Bergman fan
so that given a linear ideal $I$ with associated oriented matroid
$M_I$, the positive tropical variety associated to $I$ is equal to
the positive Bergman fan of $M_I$.

We give a proof of Sturmfels' conjecture: appropriately
subdivided, $\B^+(M)$ is a geometric realization of the order
complex of the proper part of the Las Vergnas face lattice of $M$.
$\B^+(M)$ is homeomorphic to a sphere and naturally sits
inside $\B(\underline{M})$, the Bergman complex of the underlying
unoriented matroid of $M$.  We conclude by showing that, for the
oriented matroid of the complete graph $K_n$, the face poset of a
certain ``coarse" subdivision of $\B^+(K_n)$ is dual to the face
poset of the associahedron $A_{n-2}$.

The paper is organized as follows.  In Section \ref{SectionMw} we
introduce a certain oriented matroid $M_\w$ which will play an
important role in our work.  In Section \ref{PosBerg} we define
the positive Bergman complex and prove our main theorem. In
Section \ref{PosTrop} we explain the relation between the positive
Bergman complex of an oriented matroid and the positive tropical
variety of a linear ideal. In Section \ref{Top} we describe the
topology of the positive Bergman complex of an oriented matroid.
Finally, in Sections \ref{Kn} and \ref{hooklength} 
we describe in detail the positive Bergman
complex of the oriented matroid of $K_n$: we relate it to the associahedron,
and we give a 
 formula for the number of full-dimensional 
fine cells within a full-dimensional
coarse cell.

Throughout this paper we will abuse notation and use $M$ to denote
either a matroid or oriented matroid, depending on the context.
Similarly, we will use the term ``circuits" to describe either
unsigned or signed circuits.  When the distinction between
matroids and oriented matroids is important, we will use
$\underline{M}$ to denote the underlying matroid of an oriented
matroid $M$.

\section{The Oriented Matroid $M_w$}\label{SectionMw}

Let $M$ be an oriented matroid on the ground set $[n]=\{1, 2,
\dots , n\}$ whose collection of signed circuits is $\mathscr{C}$.  Let
$\w \in \R^n$ and regard $\w$ as a weight function on $[n]$. For
any circuit $C \in \mathscr{C}$ define $\init_\w (C)$ to be the
$\w$-maximal subset of the circuit $C$ -- in other words, the
collection of elements of $C$ which have the largest weight.  We
will say that the circuit $C$ {\it achieves its largest weight
with respect to $\w$} at $\init_\w (C)$. Define $\init_\w (\mathscr{C})$
to be the collection of inclusion-minimal sets of the collection
$\{\init_\w (C) \mid C \in \mathscr{C}\}$. We then define $M_\w$ to be the
oriented matroid on $[n]$ whose collection of circuits is
$\init_\w (\mathscr{C})$.

It is not clear that $M_\w$ is a well-defined oriented matroid; we
will prove this shortly.

Given $\w\in \R^n$, let $\FF(\w)$ denote the unique flag of
subsets $\emptyset = F_0 \subset F_1 \subset \dots \subset F_k
\subset F_{k+1}=[n]$ such that $\w$ is constant on each set $F_i \setminus
F_{i-1}$ and satisfies $\w |_{F_i \setminus F_{i-1}} <\w
|_{F_{i+1}\setminus F_{i}}$. We call $\FF(\w)$ the {\it flag} of $\w$, and
we say that the {\it weight class} of $\w$ or of the flag $\FF$ is
the set of vectors $\nu$ such that $\FF(\nu) = \FF$.

It is clear that $M_\w$ depends only on the flag $\FF:=\FF(\w)$
and so we also refer to this oriented matroid as $M_\FF$.

\begin{figure}[h] 
\centering
\includegraphics[height=1in]{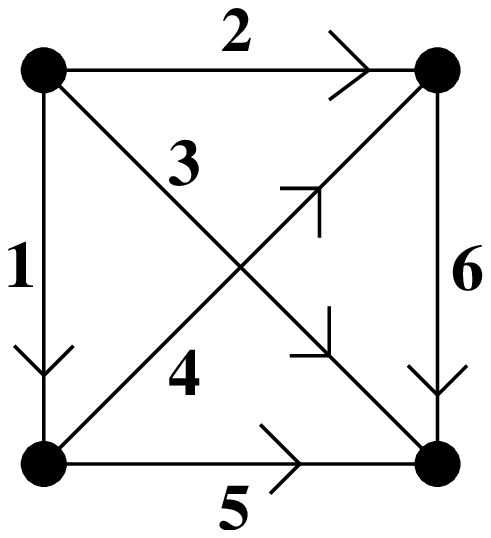}
\caption{The digraph $D$}
\label{D}
\end{figure}

\begin{figure}[h] 
\centering
\includegraphics[height=1in]{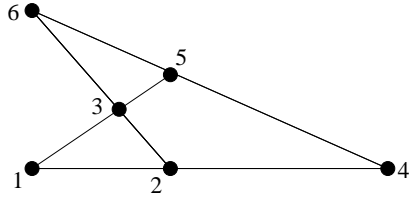}
\caption{A point configuration}
\label{Point}
\end{figure}

\begin{example}\label{MainExample}
Let $M$ be the oriented matroid of the digraph $D$ shown in
Figure~\ref{D}. Equivalently, let $M$ be the oriented matroid of the
point configuration shown in Figure~\ref{Point}.
\end{example}

%

Note that $D$ is an acyclic orientation of $K_4$, the complete
graph on $4$ vertices. The signed circuits $\mathscr{C}$ of $M$ are
$\{1\bar{2}4, 1\bar{3}5, 2\bar{3}6, 4\bar{5}6, 1\bar{2}5\bar{6},
1\bar{3}46, 2\bar{3}\bar{4}5\}$, together with the negatives of
every set of this collection. Choose $\w$ such that $\w_6 < \w_1 =
\w_2 = \w_3 = \w_4 = \w_5$, which corresponds to the flag
$\emptyset \subset \{6\} \subset \{1, 2, 3, 4, 5, 6\} $.

\begin{figure}[h] 
\centering
\includegraphics[height=1in]{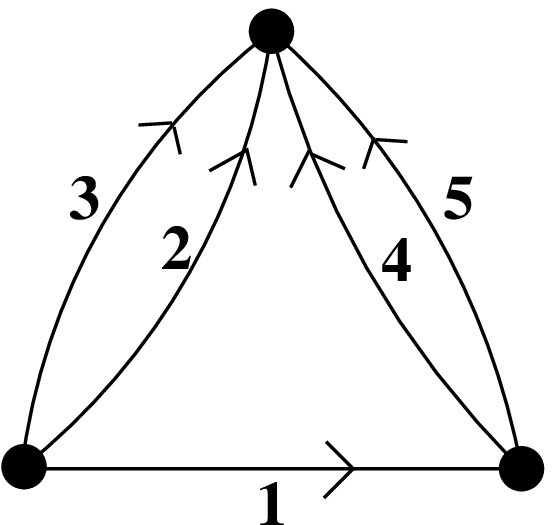}
\caption{The digraph $D'$.}
\label{Mw}
\end{figure}

If we calculate $\init_\w (C)$ for each $C \in \mathscr{C}$ we get
$\{1\bar{2}4, 1\bar{3}5, 2\bar{3}, 4\bar{5}, 1\bar{2}5, 1\bar{3}4,
2\bar{3}\bar{4}5\}$, together with the negatives of every set of
this collection.  However, $2\bar{3}\bar{4}5$ is not
inclusion-minimal in this collection, as it contains $2\bar{3}$
and $\bar{4}5$. Thus $\init_\w (\mathscr{C})$ is equal to $\{1\bar{2}4,
1\bar{3}5, 2\bar{3}, 4\bar{5}, 1\bar{2}5, 1\bar{3}4\}$, together
with the negatives of every set, and $M_\w$ is the oriented
matroid whose collection of signed circuits is $\init_\w (\mathscr{C})$.
Notice that in this case $M_\w$ is the oriented matroid of the
digraph $D^{\prime}$ in Figure \ref{Mw}.

We must show that $M_\w$ is well-defined. For convenience, we
review here the {\it circuit axioms for oriented matroids}
~\cite{RedBook}:

\renewcommand{\labelenumi}{C\theenumi.}
\begin{enumerate}
\item  $\emptyset$ is not a signed circuit. \label{C1} \item  If
$X$ is a signed circuit, then so is $-X$. \label{C2} \item  No
proper subset of a circuit is a circuit. \label{C3} \item  If
$X_0$ and $X_1$ are circuits with $X_1 \neq -X_0$ and $e \in X_0^+
\cap X_1^-$, then there is a third circuit $X\in \mathscr{C}$ with $X^+
\subseteq (X_0^+ \cup X_1^+) \setminus \{e\}$ and $X^- \subseteq
(X_0^- \cup X_1^-) \setminus \{e\}$. \label{C4}
\end{enumerate}

We will also need the following stronger characterization of
oriented matroids:

\begin{theorem}\cite[Theorem~3.2.5]{RedBook}
Let $\mathscr{C}$ be a collection of signed subsets of a set $E$
satisfying {\rm C1, C2, C3}.  Then {\rm C4} is equivalent to {\rm
C4'}:
\begin{enumerate}
\item[C4'.] for all $X_0, X_1 \in \mathscr{C}$, $e \in X_0^+ \cap X_1^-$
and $f \in (X_0^+ \setminus X_1^-) \cup (X_0^- \setminus X_1^+)$,
there is a $Z\in \mathscr{C}$ such that $Z^+ \subseteq (X_0^+ \cup X_1^+)
\setminus \{e\}$, $Z^- \subseteq (X_0^- \cup X_1^-) \setminus
\{e\}$, and $f \in Z$. \label{C4'}
\end{enumerate}
\end{theorem}

\begin{proposition}
Let $M$ be an oriented matroid on $[n]$ and $\w\in \R^n$.  Then
$M_\w$ is an oriented matroid.
\end{proposition}

\begin{proof}
The strategy of our proof is to show that if {\rm C1, C2, C3} and
{\rm C4'} hold for $M$, then {\rm C1, C2, C3} and {\rm C4} hold
for $M_\w$.   First note that it is obvious that {\rm C1} and {\rm
C2} hold for $M_\w$.  {\rm C3} holds for $M_\w$ because we defined
$\init_\w (\mathscr{C})$ to consist of {\it inclusion-minimal} elements.
It remains to show that {\rm C4} holds for $M_\w$.  To do this, we
start with two circuits in $M_\w$, lift them to circuits in $M$,
and then use {\rm C4'} for $M$ to show that {\rm C4} holds for
$M_\w$.

Take $Y_0$ and $Y_1$ in $\init_\w (\mathscr{C})$ such that $Y_1 \neq -Y_0$
and $e \in Y_0^+ \cap Y_1^-$.  By definition, there exist circuits
$X_0$ and $X_1$ of $M$ such that $Y_0 = \init_\w (X_0)$ and $Y_1 =
\init_\w (X_1)$. Notice that the presence of $e$ in $Y_0$ and
$Y_1$ guarantees that the maximum weights occurring in $X_0$ and
in $X_1$ are both equal to $\w_e$.

Choose any $f \in (Y_0^+ \setminus Y_1^-) \cup (Y_0^- \setminus
Y_1^+)$. Clearly such an $f$ exists. Then $f \in (X_0^+ \setminus
X_1^-) \cup (X_0^- \setminus X_1^+)$. By {\rm C4'} for $M$,
there exists a circuit $X$ in $\mathscr{C}$ such that
\begin{enumerate}
\item[$(a_1)$] \quad $X^+ \subseteq (X_0^+ \cup X_1^+) \setminus
\{e\}$, \item[$(a_2)$] \quad $X^- \subseteq (X_0^- \cup X_1^-)
\setminus \{e\}$, and \item[$(a_3)$] \quad $f\in X$.
\end{enumerate}

Look at $\init_\w (X)$.  We will prove that $\init_\w (X)$
contains the third circuit of $M_{\w}$ which we are looking for.
We want to show that
\begin{enumerate}
\item[$(b_1)$] \quad $\init_\w (X)^+ \subseteq (Y_0^+ \cup
                               Y_1^+) \setminus \{e\}$
\item[$(b_2)$] \quad $\init_\w (X)^- \subseteq (Y_0^- \cup
                               Y_1^-) \setminus \{e\}$
\end{enumerate}

First, it is obvious that $e$ is not in $\init_\w (X)^+$, since
$e$ was not in $X^+$.  Clearly $\init_\w (X)^+$ is a subset of
$X_0^+ \cup X_1^+$.  To show $(b_1)$ and $(b_2)$, we just need to
show that the maximum weight which occurs in $X$ is also equal to
$\w_e$.
%
By $(a_1)$ and $(a_2)$, this maximum weight
is at most $\w_e$. By $(a_3)$, equality is attained for $f \in X$,
since $\w_f=\w_e$.
%

Note that if $\init_\w (X)$ is not inclusion-minimal in the set
$\{\init_\w (C) \mid C \text{ a circuit of }M\}$, then it contains
some inclusion-minimal $\init_\w (W)$ for another circuit $W$ of
$M$. And since $\init_\w (X)^+ \subseteq (Y_0^+ \cup
                               Y_1^+) \setminus \{e\}$ and
$\init_\w (X)^- \subseteq (Y_0^- \cup
                               Y_1^-) \setminus \{e\}$, it is
clear that we also have $\init_\w (W)^+ \subseteq (Y_0^+ \cup
                               Y_1^+) \setminus \{e\}$ and
$\init_\w (W)^- \subseteq (Y_0^- \cup
                               Y_1^-) \setminus \{e\}$.

\end{proof}

In the following proposition, we describe the bases of $M_\w$ and
their orientations. Let us say that if $S \subset [n]$ and $\w \in
\R^n$, the {\it $\w$-weight} of $S$ is the sum $\sum_{i\in S}
\w_i$.

\begin{proposition}
The bases of $M_\w$ are the bases of $M$ which have minimal
$\w$-weight. The basis orientations of $M_{\w}$ are equal to their
orientations in $M$.
\end{proposition}

\begin{proof}
We know that if $N$ is an oriented matroid on $[n]$ with signed
circuits $\mathscr{C}$, then the bases of $N$ are the maximal subsets of
$[n]$ which contain no circuit.  Thus, the bases $\B$ of $M_\w$
are the maximal subsets of $[n]$ which do not contain a set of the
collection $\{ \init_\w (C) \mid C \text{ a circuit of }M\}$. We
want to show that $\B$ is exactly the set of bases of $M$ which
have minimal $\w$-weight.

First let us choose a basis $B$ of $M$ which has minimal
$\w$-weight. We claim that $B$ is independent in $M_\w$. Suppose
that $B$ contains a subset of the form $\init_\w (C)$, where $C$
is a circuit of $M$. Write $B = \init_\w (C) \cup \{b_1, \dots ,
b_m\}$, and choose any $c\in \init_\w (C)$. We now construct a new
basis $B^{\prime \prime}$ of $M$ of smaller weight than $B$, as
follows. Start with the set $C$.  We know that $C \cup \{b_1,
\dots , b_m\}$ is a spanning set of $M$.  Let $B_0$ be a minimal
subset (possibly empty) of $\{b_1, \dots , b_m\}$ such that $C
\cup B_0$ is still a spanning set of $M$. Since $C$ is minimally
dependent in $M$, the set
$B^{\prime}:=(C\setminus \{c\})\cup B_0$ will be a basis of $M$.
Now by the basis exchange axiom, for some $b\in B^{\prime}
\setminus B$, $B^{\prime \prime} := (B \setminus \{c\}) \cup
\{b\}$ is a basis of $M$.  Since $b\in B^{\prime} \setminus B$, it
follows that $b\in C\setminus \init_\w (C)$.  Thus, the weight of
$b$ is strictly less than the weight of $c$, which implies that
the weight of $B^{\prime \prime}$ is less than the weight of $B$.
This is a contradiction.

The previous argument shows that $r(M) \leq r(M_{\w})$. But every
circuit of $M_{\w}$ is contained in a circuit of $M$, so $r(M)
\geq r(M_{\w})$. It follows that $M$ and $M_{\w}$ have the same
rank.

Now let $B$ be a basis of $M_{\w}$, \emph{i.e.} $B$ is a maximal subset
of $[n]$ which does not contain an element of $\{\init_\w (C) \mid
C \text{ a circuit of }M\}$. We claim that $B$ is a basis of $M$
with minimal $\w$-weight. First note that $B$ is clearly
independent in $M$: if it were dependent in $M$, it would contain
some circuit $C$ of $M$ and hence would contain $\init_\w (C)$.
Also, $B$ has $r(M_{\w}) = r(M)$ elements. Therefore it is a basis
of $M$.

Finally, let us show that $B$ has minimal $\w$-weight.  Suppose
not. Let $c_1, \dots , c_r$ be the elements of $B$ with highest
weight. We claim that $\{c_1, \dots , c_r\} \supseteq \init_\w
(C)$ for some circuit $C$ of $M$, which will be a contradiction.
Since $B$ is a basis of $M$, adding to $B$ any element of
$[n]\setminus B$ creates a circuit. Since $B$ is not a basis of
minimal $\w$-weight, there must be an element $b\in [n]\setminus
B$ such that the weight of $b$ is strictly less than the weight of
each of the elements $c_1, \dots , c_r$.  Thus $B \cup b$ contains
a circuit $C$, and $\init_\w (C) \subseteq \{c_1, \dots , c_r\}$,
as claimed.

\medskip

To prove the claim about orientations, start with a basis $B$ of
minimal weight of $M$. Recall that an oriented matroid has exactly
two basis orientations, which are opposite to each other.
Therefore we can assume without loss of generality that $B$ has
the same orientation in $M$ and $M_{\w}$.

For any two ordered bases $B_1 = (e, x_2, \ldots, x_r)$ and $B_2 =
(f, x_2 \ldots, x_r)$ of $M_{\w}$ with $e \neq f$, we have
\[
\chi_{\w}(e, x_2, \ldots, x_r) = -C_{\w}(e)C_{\w}(f) \chi_{\w}(f,
x_2 \ldots, x_r),
\]
where $\chi_{\w}$ is the chirotope of $M_{\w}$, and $C_{\w}$ is
one of the two opposite signed circuits of $M_{\w}$ in $\{e, f,
x_2, \ldots, x_r\}$. Now $B_1$ and $B_2$ are also bases of $M$;
let $C$ be one of the two opposite signed circuits of $M$ in $\{e,
f, x_2, \ldots, x_r\}$. Then $\init_{\w}(C)$ contains a circuit of
$M_{\w}$; it must be either $C_{\w}$ or $-C_{\w}$. In any case, we
have $C(e)C(f) = C_{\w}(e)C_{\w}(f)$, so $\chi_{\w}(B_1)
\chi_{\w}(B_2)=\chi(B_1)\chi(B_2)$. It follows that if $B_1$ has
the same orientation in $M$ and $M_{\w}$, then so does $B_2$.

Recall that one can obtain any basis of a matroid from any
other by a sequence of simple basis exchanges of the type above.
Since $B$ has the same orientation in $M$ and $M_{\w}$, so does
any other basis of $M_{\w}$.
\end{proof}

\section{The Positive Bergman Complex}\label{PosBerg}

Our goal in this section is to define the positive Bergman complex
of an oriented matroid $M$ and to relate it to the Las Vergnas
face lattice of $M$, thus answering Sturmfels' question
~\cite{Sturmfels2}. We begin by giving some background on the
Bergman complex and fan of a (unoriented) matroid.

The \emph{Bergman fan} of a matroid $M$ on the ground set $[n]$ is
the set
\[
\widetilde{\B}(M) := \{ \w \in \R^n \,\, : \,\, M_{\w} \,\,
\text{has no loops} \}.
\]
The \emph{Bergman complex} of $M$ is
\[
\B(M) := \{ \w \in S^{n-2} \,\, : \,\, M_{\w} \,\, \text{has no
loops} \},
\]
where $S^{n-2}$ is the sphere $ \{\,\w \in \R^n \,\, : \,\, \w_1 +
\cdots + \w_n = 0 \, , \, \w_1^2 + \cdots + \w_n^2 = 1 \}.$

For simplicity, in this section we will concentrate on the Bergman
complex of $M$, but similar arguments hold for the Bergman fan of
$M$.

Since the matroid $M_{\w}$ depends only on the weight class that
$\w$ is in, the Bergman complex of $M$ is a disjoint union of the
weight classes of flags $\FF$ such that $M_{\FF}$ has no loops.
We say that the weight class of a flag $\FF$ is \emph{valid} for $M$ if
$M_{\FF}$ has no loops.

There are two polyhedral subdivisions of $\B(M)$, one of which is
clearly finer than the other.
The \emph{fine subdivision} of $\B(M)$ is the subdivision of
$\B(M)$ into valid weight classes: two vectors $u$ and $v$ of
$\B(M)$ are in the same class if and only if $\FF(u) = \FF(v)$.
The \emph{coarse subdivision} of $\B(M)$ is the subdivision of
$\B(M)$ into $M_{\w}$-equivalence classes: two vectors $u$ and $v$
of $\B(M)$ are in the same class if and only if $M_u = M_v$.


The following results give alternative descriptions of $\B(M)$:

\begin{theorem}\cite{CarlyFed}\label{CF}
Given an (unoriented) matroid $M$ on the ground set $[n]$ and $\w
\in \R^n$ which corresponds to a flag $\FF:=\FF(\w)$, the
following are equivalent:
\renewcommand{\labelenumi}{\theenumi.}
\begin{enumerate}
\item  $M_\FF$ has no loops. \item For each circuit $C$ of $M$,
$\init_{\w}(C)$ contains at least two elements of $C$. \item $\FF$
is a flag of flats of $M$.
\end{enumerate}
\end{theorem}

\begin{corollary}\cite{CarlyFed}\label{CF2}
Let $M$ be a (unoriented) matroid.  Then the fine subdivision of
the Bergman complex $\B(M)$ is a geometric realization of
$\Delta(L_M-\{\,\hat{0}\,,\hat{1}\,\}\,)$, the order complex of
the proper part of the lattice of flats of $M$.
\end{corollary}

We are now ready for the positive analogues of these concepts.
The \emph{positive Bergman fan} of an oriented matroid $M$ on the 
ground set $[n]$ is
\[
\widetilde{\B}^+(M) := \{ \w \in \R^n \,\, : \,\, M_{\w} \,\,
\text{is acyclic} \}.
\]
The \emph{positive Bergman complex} of $M$ is
\[
\B^+(M) := \{ \w \in S^{n-2} \,\, : \,\, M_{\w} \,\, \text{is
acyclic} \}.
\]

Within each equivalence class of the coarse subdivision of
$\B(M)$, the vectors $\w$ give rise to the same unoriented $M_{\w}$. Since
the orientation of $M_{\w}$ is inherited from that of $M$, they
also give rise to the same oriented matroid $M_{\w}$. Therefore each
coarse cell of $B(M)$ is either completely contained in or
disjoint from $\B^+(M)$. Thus $\B^+(M)$ inherits the coarse and
the fine subdivisions from $\B(M)$, and each subdivision of
$\B^+(M)$ is a subcomplex of the corresponding subdivision of
$\B(M)$.

Let $M$ be an acyclic oriented matroid on the ground set $[n]$. We
say that a covector $v \in \{+, -, 0\}^n$ of $M$ is
\emph{positive} if each of its entries is $+$ or $0$.  We say that
a flat of $M$ is \emph{positive} if it is the $0$-set of a
positive covector. Additionally, we consider the set $[n]$ to be a
positive flat.
For example, if $M$ is the matroid of Example \ref{MainExample},
then $16$ is a positive flat which is the $0$-set of the
positive covector $(0++++\,0)$.  

The {\it Las Vergnas face lattice} $\FF_{\ell v}(M)$ is the
lattice of positive flats of $M$, ordered by containment.
Note that the lattice of positive flats of the oriented matroid
$M$ sits inside $L_M$, the lattice of flats of
$M$.

\begin{example}
Let $M$ be the oriented matroid from Example \ref{MainExample}.
The positive covectors of $M$ are $\{0++++\,0, 0\,0+0++, +++\,0\,0\,0,
0+++++, +++\,0++, +++++\,0, ++++++\}$ and the positive flats are
$\{16, 124, 456, 1, 4, 6, \emptyset, 123456\}$. The lattice of
positive flats of $M$ is shown in Figure \ref{PosFlats}, alongside
the lattice of flats of $M$.

\begin{figure}[h]
\centering
\includegraphics[height=2in ]{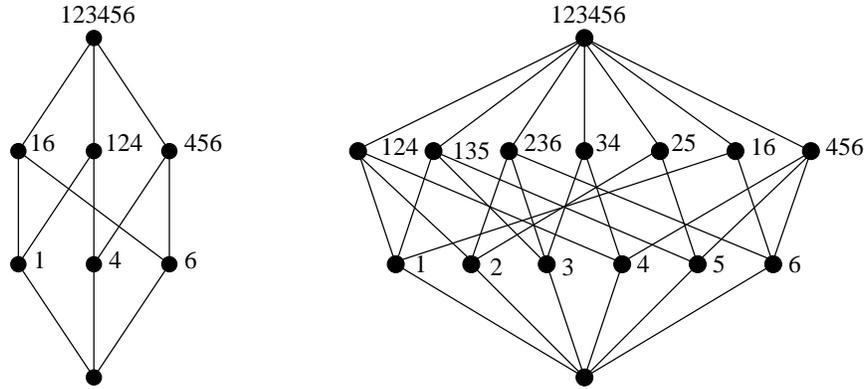}
\caption{The lattice of positive flats and the lattice of flats.}
\label{PosFlats}
\end{figure}
\end{example}

We now give an analogue of Theorem \ref{CF}.

\begin{theorem}\label{maintheorem}
Given an oriented matroid $M$ and $\w \in \R^n$ which corresponds
to a flag $\FF:=\FF(\w)$, the following are equivalent:
\renewcommand{\labelenumi}{\theenumi.}
\begin{enumerate}
\item  $M_\FF$ is acyclic. \label{1} \item For each signed circuit
$C$ of $M$, $\init_{\w}(C)$ contains a positive element and a
negative element of $C$. \label{2} \item $\FF$ is a flag of
positive flats of $M$. \label{3}
\end{enumerate}
\end{theorem}

\begin{proof}
First we will show that \ref{1} and \ref{2} are equivalent.  The
statement that $M_\w$ is acyclic means that $M_\w$ has no
all-positive circuit: in other words, each circuit of $M_\w$
contains a positive and a negative term.  Since $M_\w$ is the
matroid whose circuits are the inclusion-minimal elements of the
set $\{\init_\w (C) \mid C \text{ a circuit of }M\}$, this means
that for each circuit $C$ of $M$, $\init_\w C$ contains a positive
and a negative term.  Finally, this is equivalent to the statement
that for each circuit $C$ of $M$, $C$ achieves its maximum value
with respect to $\w$ on both $C^+$ and $C^-$.

Next we show that \ref{3} implies \ref{2}. Assume we have an $\w$
such that $\FF$ is a flag of positive flats.  Let the flats of
this flag be $F_1 \subset F_2 \subset \ldots \subset F_k$. For
each $F_i$, $([n] - F_i)$ is a positive covector.  By
orthogonality of circuits and covectors, we know that for any
circuit $C$ and any covector $Y$, $(C^+ \cap Y)$ and
$(C^- \cap Y)$ are either both empty or both non-empty. For any
circuit $C$ of $M$, consider the largest $i$ such that $C \cap
([n] - F_i)$ is non-empty.  Then clearly $C$ will attain its
maximum on $F_{i+1} - F_{i}$ and $\init_\w (C)$ contains a
positive element and a negative element of $C$.

Finally, assume that \ref{1} and \ref{2} hold, but \ref{3} does
not. From \ref{1} we know that $M_{\FF}$ is acyclic; therefore the
unoriented $M_{\FF}$ has no loops, and $\FF$ is a flag of flats by
Theorem \ref{CF}. Let $F_i$ be a flat which is not positive; 
by \cite[Proposition~9.1.2]{RedBook} this
is equivalent to saying that $M/F_i$ is not acyclic. Let $C$ be a
positive circuit of $M/F_i$; then we can find a circuit $X$ of $M$
such that $C = X - F_i$. Then $X$ has positive elements of weight
greater than $\w_i$, and no negative elements of weight greater
than $\w_i$. It follows that $\init_{\w} (X)$ is positive,
contradicting \ref{2}.

%
\end{proof}

\begin{corollary}\label{LasVergnas}
Let $M$ be an oriented matroid. Then the fine subdivision of
$\B^+(M)$ is a geometric realization of $\Delta(\FF_{\ell v}(M)
-\{\,\hat{0}\,,\hat{1}\,\}\,)$, the order complex of the proper
part of the Las Vergnas face lattice of $M$.
\end{corollary}

\section{Connection with Positive Tropical Varieties}\label{PosTrop}

In~\cite{SpeyerWilliams}, the notion of the positive part of the
tropicalization of an affine variety (or positive tropical
variety, for short) was introduced, an object which has the
structure of a polyhedral fan in $\R^n$.  In order to describe
this object, we must define an {\it initial
ideal}.

Let $\RR = \R[x_1, \ldots, x_n]$ and $\w \in \R^n$. If $f = \sum
c_i \bf{x^{a_i}} \in \RR$, define the \emph{initial form}
$\init_{\w} (f) \in \RR$ to be the sum of all terms $c_i
\bf{x^{a_i}}$ such that the inner product $\w \cdot \bf{a_i}$ is
maximal. For an ideal $I$ of $\R[x_1, \ldots, x_n]$, define the
\emph{initial ideal} $\init_\w (I)$ to be the ideal generated by
$\init_\w (f)$ for all $f \in I$.

If $I$ is an ideal in a polynomial ring with $n$ variables, the
positive tropical variety associated to $I$ is denoted by $\Trop^+
V(I)$ and can be characterized as follows:
$$\Trop^+ V(I) = \{\w\in \R^n \mid \text{
$\init_\w (I)$ contains no nonzero polynomials in $\R^+ [x_1,\dots
, x_n]$}\}.$$

Now recall that if $I$ is a linear ideal (an ideal generated
by linear forms), we can associate to it
an oriented matroid $M_I$ as follows. Write each linear form $f\in
I$ in the form $a_1 x_{i_1}+a_2 x_{i_2}+\dots + a_m x_{i_m} = b_1
x_{j_1}+b_2 x_{j_2}+\dots + b_n x_{j_n}$, where $a_i, b_i > 0$ for
all $i$. We then define $M_I$ to be the oriented matroid whose set
of signed circuits consists of all minimal collections of the form
$\{i_1 i_2 \dots i_m \bar{j}_1 \bar{j}_2 \dots \bar{j}_n \}$. We
now prove the following easy statement.

\begin{proposition}
If $I$ is a linear ideal and $M_I$ is the associated oriented
matroid, then $\Trop^+ V(I) = \widetilde{B}^+(M_I)$.
\end{proposition}
\begin{proof}
\begin{eqnarray*}
\widetilde{B}^+(M_I) &=& \{\w \in \R^n \mid (M_I)_\w \text{ is acyclic}\}\\
&=& \{\w \in \R^n \mid (M_I)_\w \text{ has no all-positive circuit}\}\\
&=& \{\w \in \R^n \mid \init_\w (I) \text{ contains no nonzero
polynomial
in }\R^+[x_1, \dots , x_n]\} \\
&=& \Trop^+ V(I).
\end{eqnarray*}
\end{proof}

\section{Topology of the Positive Bergman Complex}\label{Top}

The topology of the positive Bergman complex of an oriented
matroid is very simple: it is homeomorphic to a sphere.  This follows from
Corollary \ref{LasVergnas} together with results about the Las
Vergnas face lattice, which we will review here.

\begin{theorem}\cite[Theorem~4.3.5]{RedBook} \label{LV}
Let $M$ be an acyclic oriented matroid of rank $r$.  Then the Las
Vergnas lattice $\FF_{\ell v}(M)$ is isomorphic to the face
lattice of a $PL$ regular cell decomposition of the
$(r-2)$-sphere.
\end{theorem}

\begin{proposition}\cite[Proposition~4.7.8]{RedBook} \label{Cell}
Let $\Delta$ be a regular cell complex.  Then its geometric
realization is homeomorphic to the geometric realization of the
order complex of its face poset.
\end{proposition}

The previous two results imply that the geometric realization of
the order complex of the Las Vergnas lattice is homeomorphic to a
sphere.  

Putting this together with Corollary \ref{LasVergnas}, 
we get the following result.

\begin{corollary}
The positive Bergman complex of an oriented matroid is homeomorphic to 
a sphere.
\end{corollary}

\section{The positive Bergman complex of the complete graph}\label{Kn}

In this section, we wish to describe the positive Bergman complex
$\B^+(K_n)$ of the graphical oriented matroid $M(K_n)$ of an
acyclic orientation of the complete graph $K_n$. We start by reviewing the
description of the Bergman complex $\B(K_n)$ of the unoriented
matroid $M(K_n)$, obtained in~\cite{CarlyFed}. For the moment we
need to consider $K_n$ as an unoriented graph.

\begin{figure}[h]
\centering
\includegraphics[height=1in ]{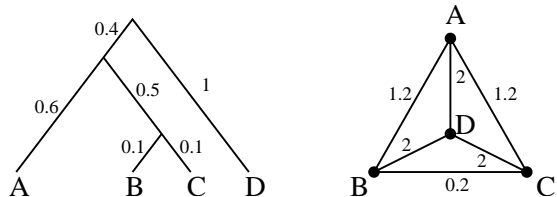}
\caption{An equidistant tree and its corresponding distance
vector.}
 \label{figure:tree}
\end{figure}

An \emph{equidistant $n$-tree} $T$ is a rooted tree with $n$
leaves labeled $1, \ldots, n$, and lengths assigned to each edge
in such a way that the total distances from the root to each
leaf are all equal. The internal edges are required to have 
positive lengths.
Figure \ref{figure:tree} shows an example of an equidistant
$4$-tree.

To each equidistant $n$-tree $T$ we assign a distance vector
$d_T \in \R^{n \choose 2}$: the distance $d_{ij}$ is equal to the
length of the path joining leaves $i$ and $j$ in $T$. Figure
\ref{figure:tree} also shows the distance vector of the tree,
regarded as a weight function on the edges of $K_4$.

The Bergman fan $\widetilde{\B}(K_n)$ can be regarded as a space
of equidistant $n$-trees, as the following theorem shows.

\begin{theorem}~\cite{CarlyFed, Semple}
The distance vector of an equidistant $n$-tree, when regarded as
a weight function on the edges of $K_n$, is in the Bergman fan
$\widetilde{\B}(K_n)$. Conversely, any point in
$\widetilde{\B}(K_n)$ is the distance vector of a unique
equidistant $n$-tree.
%
%
%
%
\end{theorem}

As mentioned earlier, the fine subdivision of $\widetilde{\B}(M)$
is well understood for any matroid $M$. The following theorem
shows that the coarse subdivision of $\widetilde{\B}(K_n)$ also
has a nice description: it is a geometric realization of the
well-studied simplicial complex of trees $T_n$, sometimes called
the Whitehouse complex \cite{Billera, Robinson}.

\begin{theorem}~\cite{CarlyFed}
Let $\w, \w' \in \widetilde{\B}(K_n)$. Let $T$ and $T'$ be the
corresponding equidistant $n$-trees. The following are equivalent:

\begin{enumerate}
\item[1.] $\w$ and $\w'$ are in the same cell of the coarse
subdivision.

%
%

\item[2.] $T$ and $T'$ have the same combinatorial type.
\end{enumerate}
\end{theorem}

\medskip

Now we return to the setting of oriented matroids. The positive
Bergman complex $\B^+(K_n)$ is defined in terms of an acyclic
orientation of $K_n$. This graph has $n!$ acyclic orientations,
corresponding to the $n!$ permutations of $[n]$. The orientation
corresponding to the permutation $\pi$ is given by $\pi_i
\rightarrow \pi_j$ for $i<j$. Clearly the $n!$ orientations of
$K_n$ will give rise to positive Bergman complexes which are equal
up to relabeling.
Therefore, throughout this section, the edges of $K_n$ will be
oriented $i \rightarrow j$ for $i<j$.

As we go around a cycle $C$ of $K_n$, $C^+$ is the set of edges
which are crossed in the forward direction, and $C^-$ is the set
of edges which are crossed in the backward direction.

\begin{proposition}
Let $\w$ be a weight vector on the edges of the oriented complete
graph $K_n$. Let $T$ be the corresponding equidistant tree. The
following are equivalent:
\begin{enumerate}
\item[1.] $\w$ is in $\B^+(K_n)$.

\item[2.] $T$ can be drawn in the plane without crossings in such
a way that its leaves are numbered $1,2, \ldots, n$ from left to
right.

\end{enumerate}
\end{proposition}

\begin{proof}
We add three intermediate steps to the equivalence:

\begin{enumerate}

\item[$(a)$] \emph{In any cycle $C$, the $\w$-maximum is achieved
in $C^+$ and $C^-$.}

\item[$(b)$] \emph{In any triangle $C$, the $\w$-maximum is
achieved in $C^+$ and $C^-$.}

\item[$(c)$] \emph{For any three leaves $i<j<k$ in $T$, the leaf
$j$ does not branch off before $i$ and $k$; \emph{i.e.},their
branching order is one of the following:}
\begin{figure}[h] 
\centering
\includegraphics[height=.7in]{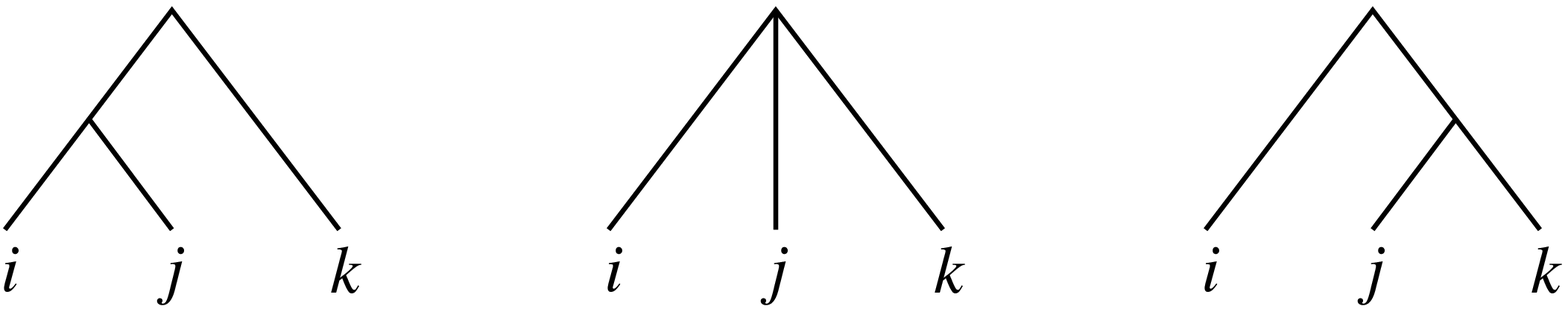}
\label{figure:branching}
\end{figure}
\end{enumerate}

The equivalence $1 \Leftrightarrow (a)$ follows from Theorem
\ref{maintheorem}, and the implication $(a) \Rightarrow (b)$ is
trivial. Now we show that $(b) \Rightarrow (a)$. Proceed by
contradiction. Consider a cycle $C=v_1 \ldots v_k$, with $k$
minimal, such that $(a)$ is not satisfied. Consider the cycles
$T=v_1v_{k-1}v_k$ and $C'=v_1v_2\ldots v_{k-1}$, which do satisfy
$(a)$. Since $C$ does not satisfy $(a)$, the edge $v_1v_{k-1}$
must be $\w$-maximum in $T$, along with another edge $e$ of the
opposite orientation. Similarly, the edge $v_{k-1}v_1$ must be
$\w$-maximum in $C'$, along with another edge $f$ of the opposite
orientation. Therefore, in $C$, the edges $e$ and $f$ are
$\w$-maximum and have opposite orientations. This is a
contradiction.

Let us now show $(b) \Leftrightarrow (c)$. In triangle $ijk$
(where we can assume $i<j<k$), $(b)$ holds if and only if we have
one of the following:
\[
\w_{ij}<\w_{jk}=\w_{ik}, \,\,\,\,\, \text{or} \,\,\,\,\,
\w_{ij}=\w_{jk}=\w_{ik}, \,\,\,\,\, \text{or} \,\,\,\,\,
\w_{jk}<\w_{ij}=\w_{ik}.
\]
These three conditions correspond, in that order, to the three
possible branching orders of $i,j$ and $k$ in $T$ prescribed by
condition $(c)$.

Finally we show $(c) \Leftrightarrow 2$. The backward implication
is immediate. We prove the forward implication by induction on
$n$. The case $n=3$ is clear. Now let $n \geq 4$, and assume that
condition $(c)$ holds. Consider a lowest internal node $v$; it is
incident to several leaves, which must have consecutive labels
$i,i+1,\ldots,j$ by $(c)$. Let $T'$ be the tree obtained from $T$
by removing leaf $i$. This smaller tree satisfies $(c)$, so it can be
drawn in the plane with the leaves in order from left to right.
Now we simply find node $v$ in this drawing, and attach leaf $i$
to it, putting it to the left of all the other leaves incident to
$v$. This is a drawing of $T$ satisfying 2.
\end{proof}

The \emph{associahedron} $A_{n-2}$ is a well-known $(n-2)$-dimensional
polytope whose vertices correspond to planar rooted trees \cite{Ziegler}.
There is a close relationship between $\B^+(K_n)$ and $A_{n-2}$.

\begin{corollary}
The face poset of the coarse subdivision of $\B^+(K_n)$, with a
$\widehat{1}$ attached, is dual to the face poset of the
associahedron $A_{n-2}$.
\end{corollary}

\begin{proof}
In the trees corresponding to the cells of $\B^+(K_n)$, the
labeling of the leaves always increases from left to right. We
can forget these labels and obtain the usual presentation of the dual to the
associahedron, whose facets correspond to planar rooted trees.
\end{proof}

\begin{figure}[h] 
\centering
\includegraphics[height=2in]{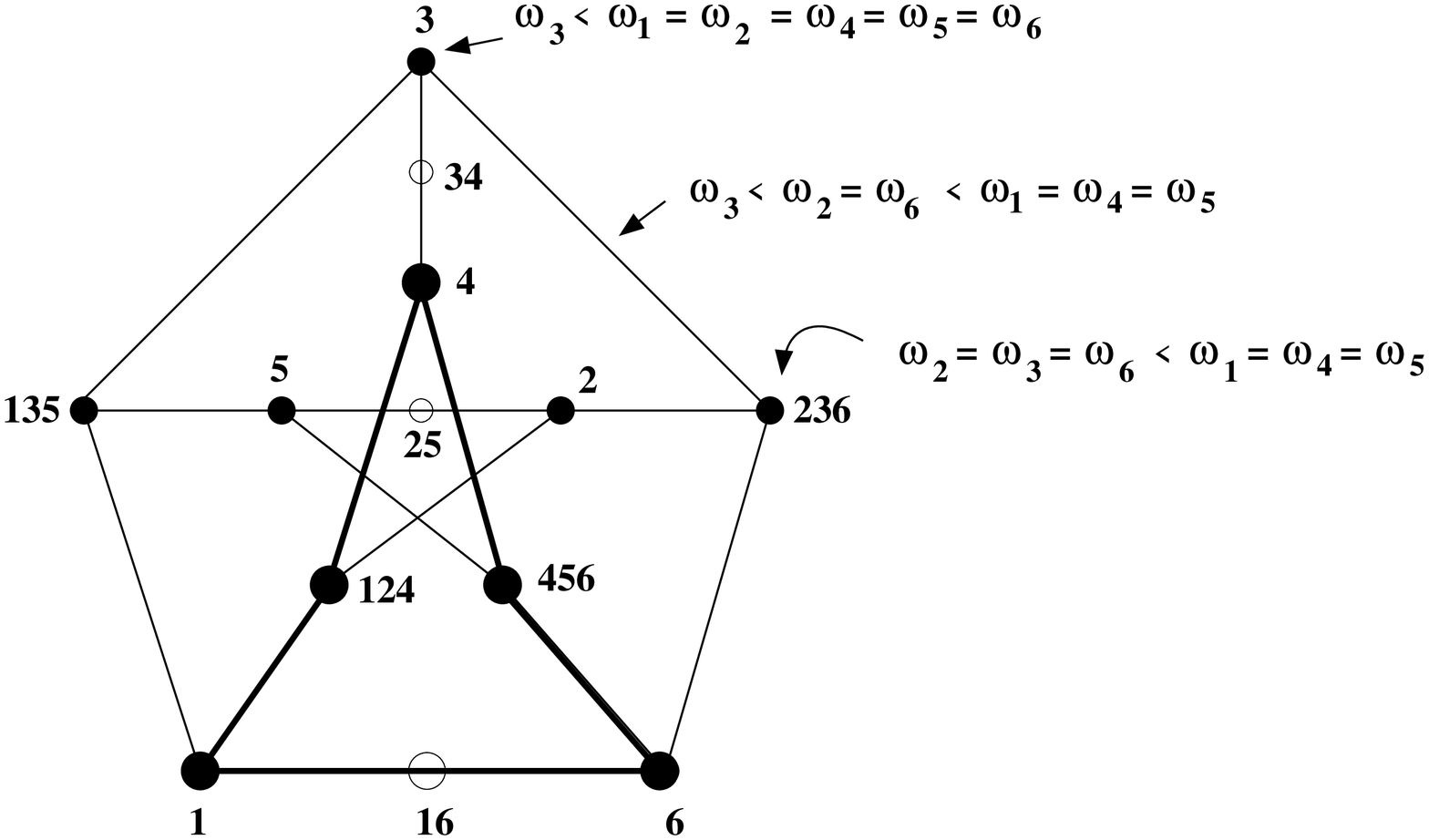}
\caption{$\B^+(K_4)\subset \B(K_4)$}
\label{Assoc}
\end{figure}

Figure \ref{Assoc} shows the positive Bergman complex of $K_4$ (in bold)
within the Bergman complex of $K_4$.  Vertices of the coarse subdivision are
shown as black circles; vertices of the fine
subdivision but not the coarse subdivision are shown as transparent
circles.  Observe that the coarse subdivision of 
$\B^+(K_4)$ is a pentagon, whose face poset
is the face poset of the associahedron $A_2$ (which is self-dual).

Now, recall that different orientations of $K_n$ give rise to
different positive Bergman complexes. Let us make two
comments about the way in which these positive Bergman complexes
fit together.

Consider the $n!$ different acyclic orientations $o(\pi)$ of
$K_n$, each corresponding to a permutation $\pi$ of $[n]$. Each
orientation $o(\pi)$ gives rise to a positive Bergman complex: it
consists of those weight vectors such that the corresponding tree
can be drawn with the leaves labeled $\pi_1, \ldots, \pi_n$ from
left to right. Clearly, each permutation and its reverse give the
same positive Bergman complex. The $\frac{n!}2$ possible positive
Bergman complexes $\B^+(K_n)$ give a covering of $\B(K_n)$, and
each one of them is dual to the associahedron $A_{n-2}$. This
corresponds precisely to the known covering of the space of trees
with $\frac{n!}2$ polytopes dual to the associahedron, as
described in \cite{Billera}.

Also, recall from \cite{CarlyFed} that the Bergman complex
$\B(K_n)$ is topologically a wedge of $(n-1)!$ spheres. It is now
easy to realize this statement combinatorially as follows. Any
tree can be drawn in the plane so that the leftmost leaf is
labeled $1$. Therefore, $\B(K_n)$ is covered by the positive
Bergman complexes corresponding to the $(n-1)!$ permutations $\pi$
such that $\pi_1=1$. It follows that $\B(K_n)$ is the wedge of
$(n-1)!$ spherical polytopes dual to the associahedron $A_{n-2}$.


\section{The number of fine cells in $\B^+(K_n)$ and
$\B(K_n)$.}\label{hooklength}

Since $\B^+(K_n)$ and $\B(K_n)$ are $(n-2)$-dimensional, we will
call the $(n-2)$-dimensional cells inside them
\emph{full-dimensional}.  In this section we will give a 
formula reminiscent of the ``hook-length'' formula for the number of full-dimensional fine cells within a
full-dimensional coarse cell of $\B(K_n)$.

\begin{proposition} \label{finecells}
Let $\tau$ be a rooted binary tree with $n$ labeled leaves. For
each internal vertex $v$ of $\tau$, let $d(v)$ be the number of
internal vertices of $\tau$ which are descendants of $v$,
including $v$. Let $C(\tau)$ be the coarse cell of $\B(K_n)$
corresponding to tree $\tau$. There are exactly
\[
\frac{(n-1)!}{\prod_v d(v)}
\]
full-dimensional fine cells in $C(\tau)$.
\end{proposition}

\begin{proof}

The cell $C(\tau)$ consists of the distance vectors $d \in \R^{n
\choose 2}$ of all equidistant $n$-trees $T$ of combinatorial type
$\tau$. Notice that $d_{ij} = 2h-2h(v)$, where $v$ is the lowest
common ancestor of leaves $i$ and $j$ in $T$, $h(v)$ is the
distance from $v$ to the root of $T$, and $h$ is the distance from
the root of $T$ to any of its leaves.

To specify a full-dimensional fine cell in $C(\tau)$, one needs to
specify the relative order of the $d_{ij}$'s. Equivalently, in the
tree $T$ that $d$ comes from, one needs to specify the relative
order of the heights of the internal vertices, consistently with
the combinatorial type of tree $\tau$. Therefore, the fine cells
in $C(\tau)$ correspond to the labellings of the $n-1$ internal
vertices of $\tau$ with the numbers $1,2,\ldots, n-1$, such that
the label of each vertex is smaller than the labels of its
offspring. In the language of \cite[Sec. 1.3]{EC1}, these are
precisely the \emph{increasing binary trees} of type $\tau'$,
where $\tau'$ is the result of removing the leaves of tree $\tau$,
and the edges incident to them. Figure \ref{figure:increasing}
shows a tree type $\tau$ and one of the increasing binary trees of
type $\tau'$.
\begin{figure}[h]
\centering
\includegraphics[height=1.5in]{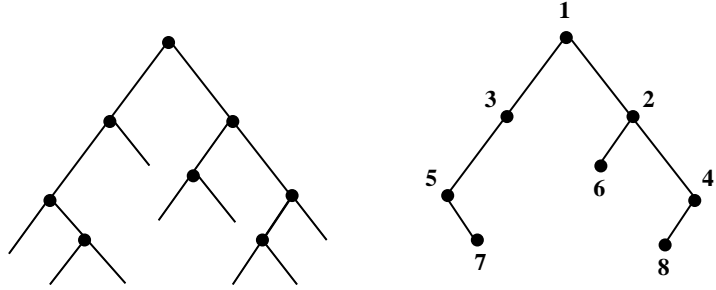}
\caption{A tree type $\tau$ and an increasing binary tree of type $\tau'$.}
\label{figure:increasing}
\end{figure}


Suppose we choose one of the $(n-1)!$ labellings of $\tau'$
uniformly at random. Let $A_{\tau'}$ be the event that the chosen
labeling $L$ is increasing; it remains to show that $P(A_{\tau'})
= 1/{\prod_{v} d(v)}$.

Let $\tau_1'$ and $\tau_2'$ be the left and right subtrees of
$\tau'$. Let $B_1$ and $B_2$ be the events that $\tau_1'$ and
$\tau_2'$ are labeled increasingly in $L$, and let $B$ be the
event that the root of $\tau$ is labeled $1$. Then $A_{\tau'} = B
\cap B_1 \cap B_2$. It is clear that $B, B_1$ and $B_2$ are
independent events. It is also clear that $P(B_1) =
P(A_{\tau_1'})$ and $P(B_2) = P(A_{\tau_2'})$. Therefore,
\begin{eqnarray*}
P(A_{\tau}) &=& P(B) P(B_1) P(B_2) \\
&=& \frac{1}{n-1} P(A_{\tau_1'}) P(A_{\tau_2'}).
\end{eqnarray*}
The result follows by induction.
\end{proof}

It is also possible to obtain analogous formulas for the number of
fine cells inside a lower-dimensional coarse cell, corresponding
to a rooted tree which is not binary. We omit the details.

Notice that Proposition \ref{finecells} is essentially
equivalent to the formula for the number of linear 
extensions of a poset whose Hasse diagram is a tree 
\cite[Supp. Ex. 3.1]{EC1}.

\begin{corollary}\label{totalfinecells}
The positive Bergman complex $\B^+(K_n)$ contains exactly $(n-1)!$
full-dimensional fine cells. The Bergman complex $\B(K_n)$
contains exactly $n!\,(n-1)!\,/\,2^{n-1}$ full-dimensional fine
cells.
\end{corollary}

\begin{proof}
We recall the known bijection between increasing binary trees with
vertices labeled $a_1< \ldots < a_k$, and permutations of $\{a_1,
\ldots, a_k\}$ \cite[Sec. 1.3]{EC1}. It is defined recursively:
the permutation $\pi(T)$ corresponding to the increasing binary
tree $T$ is $\pi(T) = \pi(T_1)\,a_1\,\pi(T_2)$, where $T_1$ and
$T_2$ are the left and right subtrees of $T$. For example, the
tree of Figure \ref{figure:increasing} corresponds to the
permutation $57316284$. It is not difficult to see how $T$ can
be recovered uniquely from $\pi(T)$.

Since the full-dimensional fine cells of $\B^+(K_n)$ are in
correspondence with the increasing binary trees with labels $1,
\ldots, n-1$, the first result follows.

\medskip

To show the second result, recall that the Bergman complex
$\B(K_n)$ is covered by $n!$ positive Bergman complexes. Each
permutation $\pi$ of $[n]$ gives rise to a positive Bergman
complex $\B^+(K_n)$; this complex parameterizes those trees which
can be drawn in the plane so that its leaves are in the order
prescribed by $\pi$. With $(n-1)!$ fine cells in each positive
Bergman complex, we get a covering of $\B(K_n)$ with $n!(n-1)!$
fine cells. Each fine cell appears several times in this covering,
since it sits inside several positive Bergman complexes.

More precisely, each binary tree with $n$ labeled leaves can be
drawn in the plane in exactly $2^{n-1}$ ways: at each internal
vertex, we may or may not switch the left and right subtrees.
Therefore, each fine cell of the Bergman complex $\B(K_n)$ is
inside $2^{n-1}$ different positive Bergman complexes. The desired
result follows.
\end{proof}

Recall that the maximum-dimensional fine cells of $\B(K_n)$ 
correspond to the maximal
chains in the lattice of flats of $K_n$; \emph{i.e.}, the partition
lattice $\Pi_n$. Thus we have given an alternative proof of the
fact that there are $n! \, (n-1)! \, / \, 2^{n-1}$ maximal chains
in $\Pi_n$ \cite[Supp. Ex. 3.3]{EC1}.

As an illustration of Corollary \ref{totalfinecells}, notice that,
in Figure \ref{Assoc}, the positive Bergman complex $\B^+(K_4)$
consists of $3! = 6$ fine cells, while the Bergman complex
$\B(K_4)$ consists of $4! \,3!\,/\,2^3 = 18$ fine cells.

\medskip

\textsc{Acknowledgments} We are grateful to Bernd Sturmfels for
suggesting this project to us, and to Lou Billera for helpful
conversations.

\raggedright


\addcontentsline{toc}{section}{References}
\begin{thebibliography}{2}

\bibitem{CarlyFed} F. Ardila and C. Klivans.
\emph{The Bergman complex of a matroid and phylogenetic trees.}
Preprint, \textsf{arXiv:math.CO/0311370}, 2003.

\bibitem{Bergman}
G. Bergman.  \emph{The logarithmic limit-set of an algebraic
variety.} Transactions of the American Mathematical Society, vol.
157, 459-469, 1971.
%

\bibitem{Billera}
L. Billera, S. Holmes, and K. Vogtmann.  \emph{Geometry of the
space of phylogenetic trees.}  Advances in Applied Mathematics,
vol. 27, 733-767, 2001.


%

\bibitem{RedBook}
A. Bj\"{o}rner, M. Las Vergnas, B. Sturmfels, N. White, G.
Ziegler. \emph{Oriented matroids.}  Encyclopedia of Mathematics
and Its Applications, vol. 46.  Cambridge University Press,
Cambridge, 1993.
%
%
%
\bibitem{Einsiedler}
M. Einsiedler, M. Kapranov, D. Lind, and T. Ward.
\emph{Non-Archimedean amoebas.} Preprint, 2003.
%

\bibitem{Lusztig1} G. Lusztig. \emph{Introduction to total positivity}, in
Positivity in Lie theory: open problems, ed. J. Hilgert, J.D.
Lawson, K.H. Neeb, E.B. Vinberg, de Gruyter Berlin, 1998, 133-145.

\bibitem{Lusztig2} G. Lusztig. \emph{Total positivity in partial flag manifolds},
Representation Theory, 2 (1998) 70-78.


%
%
\bibitem{Robinson}
A. Robinson and S. Whitehouse.  \emph{The tree representation of
$\sigma_{n+1}$.}  J. Pure Appl. Algebra, vol. 111, 245-253, 1996.

\bibitem{Semple}
C. Semple and M. Steel. \emph{Phylogenetics.} Oxford Lecture
Series in Mathematics and Its Applications, Oxford University
Press, 2003.

\bibitem{Speyer}
D. Speyer and B. Sturmfels. \emph{The tropical Grassmannian.}
Preprint, \textsf{arXiv:math.AG/0304218}, 2003.

\bibitem{SpeyerWilliams}
D. Speyer and L. Williams.\emph{The tropical totally positive
Grassmannian.} Preprint, \textsf{arXiv:math.CO/0312297}, 2003


\bibitem{EC1}
R. P. Stanley. \emph{Enumerative combinatorics, vol.1.} Cambridge
University Press, New York, 1986.

\bibitem{Sturmfels}
B. Sturmfels. \emph{Solving systems of polynomial equations.} CBMS
Regional Conference Series in Mathematics, vol. 97.  American
Mathematical Society, Providence, RI, 2002.

\bibitem{Sturmfels2}
B. Sturmfels.  personal communication, 2003.
%
%
%

\bibitem{Ziegler}
G. Ziegler. \emph{Lectures on polytopes}. Graduate Texts in Mathematics, 
vol. 152, Springer-Verlag, New York, 1995.

\end{thebibliography}
\end{document}